\documentclass[13pt,oneside]{amsart}
\usepackage{amsmath,ifthen, amsfonts, amssymb, srcltx,
   amsopn, textcomp}

\usepackage{dashrule}
\usepackage{tikz}
\usepackage{tikz-cd}
\usepackage{hyperref}
\usepackage{graphicx}

\usepackage[small]{caption}

\newcommand{\showcomments}{yes}
\renewcommand{\showcomments}{no}

\newcommand{\hidetodo}[1]
{\ifthenelse{\equal{\showcomments}{yes}}%
{#1}
}

\newsavebox{\commentbox}
\newenvironment{com}%
{\ifthenelse{\equal{\showcomments}{yes}}%
{\footnotemark
        \begin{lrbox}{\commentbox}
        \begin{minipage}[t]{1.25in}\raggedright\sffamily\tiny
        \footnotemark[\arabic{footnote}]}
{\begin{lrbox}{\commentbox}}}%
{\ifthenelse{\equal{\showcomments}{yes}}%
{\end{minipage}\end{lrbox}\marginpar{\usebox{\commentbox}}}
{\end{lrbox}}}

\newtheorem{thm}{Theorem}[section]
\newtheorem{lem}[thm]{Lemma}

\newtheorem{prop}[thm]{Proposition}

\theoremstyle{definition}
\newtheorem{defn}[thm]{Definition}
\newtheorem{rem}[thm]{Remark}

\begin{document}

\title[On unique factorization of nonperiodic words ]{On unique factorization of nonperiodic words}
\author{Brahim Abdenbi}
           \address{ }
\email{brahim.abdenbi@mail.mcgill.ca}
          \address{Dept. of Math. \& Stats.\\
                    McGill Univ. \\
                    Montreal, Quebec, Canada H3A 0B9}

\date{\today}
\thanks{Research supported by NSERC}

\maketitle

\begin{com}
{\bf \normalsize COMMENTS\\}
ARE\\
SHOWING!\\
\end{com}

\begin{abstract}
Given a bi-order $\succ$ on the free group $\mathcal{F}$, we show that every non-periodic cyclically reduced word $W\in \mathcal{F}$ admits a maximal ascent that is uniquely positioned. This provides a cyclic permutation of $W'$ that decomposes as $W'=AD$ where $A$ is the maximal ascent and $D$ is either trivial or a descent. We show that if $D$ is not uniquely positioned in $W$, then it must be an internal subword in $A$. Moreover, we show that when $\succ$ is the Magnus ordering,  $D=1_\mathcal{F}$ if and only if $W$ is monotonic. 
 \end{abstract}

\section{Introduction}
We investigate the relationship between orderability and structural properties of nonperiodic words. Weinbaum conjectured in \cite{Weinbaum90} that any nonperiodic word  $W$ of length $>1$ has a cyclic permutation that is a concatenation $UV$ where each of $U$ and $V$ appear exactly once as a prefix of a cyclic permutation of $W$ and $W^{-1}$. This conjecture was proved by Duncan-Howie in \cite{DuncanHowie92} using the right-orderability of one-relator groups \cite{BurnsHale72}. It was also proved in \cite{MR2217836} using the Critical Factorization Theorem \cite{MR501107} (also see \cite{MR1125926} and \cite{MR1958246}) and basic properties of Lyndon words \cite{MR64049}. This provided the motivation to investigate this question without the complex machinery used above, but rather using only bi-orderability of the free group. For this purpose, we introduce the notions of \textit{ascents/descents} in a cyclically reduced words. Our result provides a partial proof of Weinbaum's conjecture and gives some additional insights into the structure of cyclically reduced nonperiodic words.
The main result of this note is:

\begin{thm}\label{thm:magnus unique position} Let $X=\left\{x_1, x_2\right\}$ be an alphabet and let $W\in \mathcal{F}\left(X\right)$ be a cyclically reduced nonperiodic word of length $>1$. Then $W$ has a cyclic permutation $W'=AD$ where:
\begin{enumerate}
\item $A$ is the uniquely positioned maximal ascent in $W$.
\item If $D$ is not uniquely positioned, then it appears as an internal subword of $A$.
\item Using the Magnus ordering on $\mathcal{F}$, we have $D=1_\mathcal{F}$ if and only if $W$ is monotonic.
\end{enumerate}
\end{thm}

Theorem~\ref{thm:magnus unique position} shows that nonperiodic cyclically reduced words have cyclic permutations that factor as concatenations of maximal ascents and descents. The maximal ascents are always uniquely positioned, whereas the descents are not necessarily so. We show that when the descents are not uniquely positioned, they appear as internal subwords of the maximal ascents. Consequently, when the descents are not uniquely positioned, they have shorter lengths than maximal ascents. We also show that when the Magnus bi-ordering is used, we can assert that the ascents are equal to $W$ if and only if $W$ is monotonic, in the sense that it only contains letters in $X$ or $X^{-1}$, but not both. This is all achieved without using any of the machinery of right-orderability/local indicability of one-relator groups, but rather using only bi-orderability of the free group and some basic combinatorial arguments.   

\section{Ascents and Descents}
\begin{figure}\centering
\includegraphics[width=.7\textwidth]{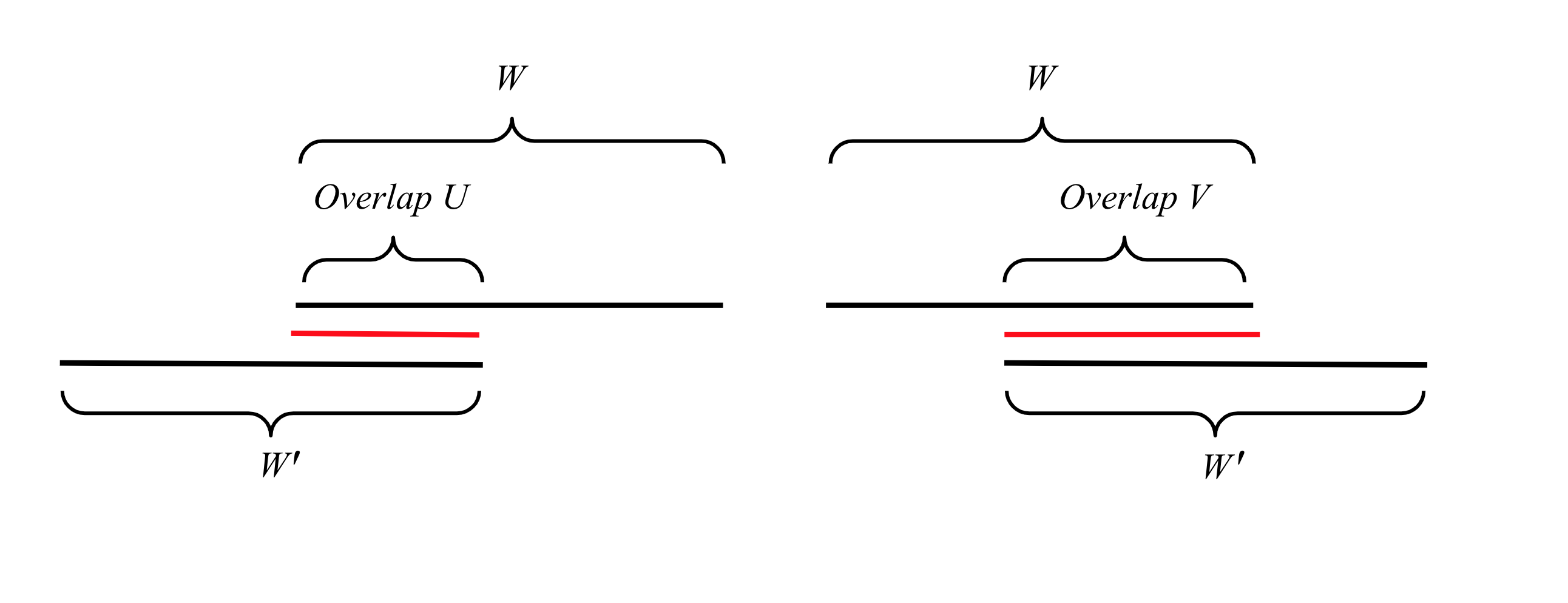}
\caption[]{\label{jump2}
We represent words as line segments. Here we see $W$ and $W'$ overlap along $U$ and $V$.}
\end{figure}

Let $X$ be an alphabet and let $W=y_1\cdots y_n$ be a word in $X^*=X\cup X^{-1}$. We say $W$ is \textit{reduced} if $y_i\neq y_{i+1}^{-1}$ for all $1\leq i<n$. It is \textit{cyclically reduced} if $W$ is reduced and $x_1\neq x_n^{-1}$. We henceforth consider only cyclically reduced words. Each $W$ represents an element $g$ in the free group $\mathcal{F}=\mathcal{F}\left(X\right)$,\textcolor{black}{and each $g\in \mathcal{F}$ is represented by a unique reduced word $W$}. For simplicity, we shall use $W$ to denote both the element of $\mathcal{F}$ and its representation in $X^*$. The \textit{spelling} of $W$ is $W=y_1\cdots y_n$ where $y_i\in X^*$. The \textit{empty} word is denoted by $1_\mathcal{F}$. The \textit{length} of $W$, denoted by $|W|$, is $n$ if $W=y_1\cdots y_n$ for $y_i\in X^*$. 
\textcolor{black}{We also represent a word $W$ as a finite graph denoted by $\overline{W}$, that is linear, directed, and labeled in $X^*$. The empty word is then represented by a single vertex.} A word $V$ is a \textit{subword} of $W$ if $W=SVU$ for some reduced words $S$ and $U$ with $|W|=|S|+|V|+|U|$.\textcolor{black}{Note that each subword $V$ of $W$ is represented by a connected subgraph $\overline{V}\subset \overline{W}$} .The subword $V$ is a \textit{prefix} if $S=1_{\mathcal{F}}$ and a \textit{suffix} if $U=1_{\mathcal{F}}$. If both $S\neq 1_{\mathcal{F}}$ and $U\neq 1_{\mathcal{F}}$, then $V$ is \textit{internal} in $W$. Subwords $U$ and $V$ of $W$ are \textit{equivalent} if they have the same spelling\textcolor{black}{and appear in different positions in $W$}. We write $U\equiv V$. Two subwords of $W$ \textit{overlap} if a nonempty suffix of one is a prefix of the other and neither is a subword of the other. We shall represent words and overlaps diagrammatically as line segments. See Figure~\ref{jump2}. 
$W$ is \textit{periodic} if there is a cyclically reduced word (\textit{period}) $U\in \mathcal{F}$ with $U^k=W$ for some $k>1$.\textcolor{black}{If $U$ is a prefix of $W$, with $W=UV$, then the word $W'=VU$ is a \textit{cyclic permutation} of $W$. Working in the free group, we have $W'=U^{-1}WU$, where $W'$ is reduced if and only if $W$ is cyclically reduced}. Let $R_{_W}$ be the set of cyclic permutations of $W$ and $W^{-1}$. A nontrivial reduced word $U\in \mathcal{F}$ is a \textit{uniquely positioned} in $W$ if 
$U$ is the prefix of exactly one element of $R_{_W}$. For example, $aa$ is uniquely positioned in $W=baaba$. The word $aba$ is not uniquely positioned in $W$ since it is a prefix in both $ababa$ and $abaab$.

A \textit{bi-ordered} group is a pair $(\mathcal{G},\prec)$ where $\mathcal{G}$ is a group and $\prec$ is a total order on $\mathcal{G}$ that is invariant under both right and left group\textcolor{black}{translation}. A well known result of Shimbireva \cite{MR0020558} states that the free group on two generators (and so every non-abelian free group) is bi-orderable. This result is sometimes attributed to Vinogradov and Magnus as well \cite{DNR}. 
\begin{defn} Let $\mathcal{F}=\mathcal{F}\left(X\right)$ be the free group on an alphabet $X$. Let $\prec$ be a bi-order on $\mathcal{F}$. A word $U\in \mathcal{F}$ is an \textit{ascent} if each prefix and each suffix of $U$ is $\succ1_\mathcal{F}$, and $U$ is a \textit{descent} if each prefix and each suffix of $U$ is $\prec1_\mathcal{F}$. The \textit{maximal ascent} $A$ of $W$ is the greatest ascent over all subwords of $R_{_W}$. 

\textcolor{black}{The} \textit{peak} (resp. \textit{low}) in $W$ is the largest (resp. smallest) prefix of $W$ with respect to $\prec$.
\end{defn}
\begin{rem}\label{rem:no overlap} Let $W$ be a cyclically reduced word representing an element in $\mathcal{F}$. Let $\prec$ be a bi-order on $\mathcal{F}$. Then:
\begin{enumerate}
\item The maximal ascent of $W$ is the largest subword in $R_{_W}$ with respect to $\prec$.
\item The inverse of an ascent is a descent.
\item If two ascents\textcolor{black}{in a given word} overlap, then the overlap is also an ascent.
\item Ascents and descents\textcolor{black}{in a given word} have no overlaps.
\item\textcolor{black}{Ascents and descents are always nontrivial. Peaks and lows, on the other hand can be trivial.}
\end{enumerate}
\end{rem}
\begin{lem}\label{lem: ascent is M minus m} Let $W$ be a  cyclically reduced word of $(\mathcal{F}(X),\prec)$. Let $A$ be the maximal ascent in $W$. Let $M$ and $m$ be the peak and low of $W$. If $A$ is a subword of $W$, then $W\succ1_\mathcal{F}$ and $M=mA$. Consequently, $W$ has exactly one subword equivalent  to $A$ and no subword of $W^{-1}$ is equivalent to $A$. Moreover, if $W=AD$ with $D\neq 1_\mathcal{F}$, then $D$ is a descent.
\end{lem}
\begin{proof}
Let $W=PAQ$, where $P$ and $Q$ are subwords of $W$. If $W\prec1_\mathcal{F}$, then $A\prec P^{-1}Q^{-1}$. But $P^{-1}Q^{-1}$ is an initial subword of a conjugate of $W^{-1}$. Indeed, we have $W^{-1}=Q^{-1}A^{-1}P^{-1}$. Then 
$AQ\left(Q^{-1}A^{-1}P^{-1}\right)Q^{-1}A^{-1}=P^{-1}Q^{-1}A^{-1}$ is a cyclic permutation of $W^{-1}$. This leads to a contradiction since $A$ is the maximal ascent in $W$. Thus, if $A$ is a subword of $W$, then $W\succ1_\mathcal{F}$ and so $A$ is not a subword of conjugates of $W^{-1}\prec1_\mathcal{F}$.

Let $W=y_1\cdots y_n$ and let $g_i=y_1\cdots y_i$. Suppose $g_j=g_iA$ for some $0\leq i, j \leq n$. Then $g_j\preceq M$ and $g_i\succeq m$. So  $g_i^{-1}\preceq m^{-1}$. Since $\prec$ is a bi-order, $A=g_i^{-1}g_j\preceq m^{-1}M$.

We now show that $m^{-1}M$ is a subword of $W$. Note that $m^{-1}M$ is a subword of $W$ whenever $|m|<|M|$. Let $m=y_1\cdots y_s$ and $M=y_1\cdots y_t$. Then $s\neq t$. Suppose $s > t$. Then $m^{-1}M=y_s^{-1}\cdots y_{t+1}^{-1}$ is a subword in $W^{-1}$. By the maximality of $A$, we have $A\preceq m^{-1}M\Rightarrow A= m^{-1}M$. So 
$A$ is a subword of $W^{-1}$,  which is a contradiction. Thus, $s<t$ and  $A=m^{-1}M$ is a subword of $W$.

It remains to show that $A$ does not appear twice in $W$. Following \cite{DuncanHowie92}, if $g_j=g_iA$, then: $$g_j\ =\ g_i\left(g_i^{-1}g_j\right)\ \succeq\ m\left(g_i^{-1}g_j\right)\ =\ m(m^{-1}M)\ =\ M$$
and so $g_j=M$ and $g_i=m$. Thus $A$ appears exactly once in $W$.

Suppose $W=AD$ where $A$ is the maximal ascent in $W$\textcolor{black}{and $D\neq 1_{\mathcal{F}}$}. Then $D\prec1_\mathcal{F}$ since otherwise $AD>A$ contradicting the maximality of $A$. If $D$ has a prefix $U\succ 1_\mathcal{F}$, then $AU\succ A$ which is a contradiction. If $D$ has a suffix $U\succ 1_\mathcal{F}$ then  $UA\succ A$ contradicting the maximality of $A$. Thus $D$ is a descent.
\end{proof}

\begin{rem}\label{rem:W square} When $W=AD$, the peak $M$\textcolor{black}{is} $A$ and the low $m$\textcolor{black}{is} $1_\mathcal{F}$.

Each cyclic permutation of $W=AD$ appears as a subword of $W^n$ with $n\geq 2$. To show that $A$ is uniquely positioned in $W$, it suffices to show that $W^2$ has exactly $2$ subwords equivalent to $A$.\textcolor{black}{That is, $A$ appears as a subword of $W^2$ in only the expected positions which are $W^2=ADAD$. In general, there are $n\geq 2$ occurrences of the subword $A$ in} $W^n=\underbrace{AD\cdots AD}_{n\ \text{times}}$.\textcolor{black}{To prove the main theorem, we will show that $W^2$ contains exactly two occurrences of the maximal ascent $A$.} 
\end{rem}
\begin{defn}\label{defn:cascade}
Let $W=AD$ and $W'=AD'$ be overlapping cyclic permutations in $W^n$ where $W'=U^{-1}WU$ and $n\geq 2$. A \textit{cascade} in $W^n$ induced by the \textit{shift} $U$ is a sequence of concatenated subwords $\left\{U_i\right\}_{i\geq 0}$ in $W$ where $U_i\equiv U$ for each $i$. 
See Figure~\ref{jump1}. 
\begin{figure}\centering
\includegraphics[width=.7\textwidth]{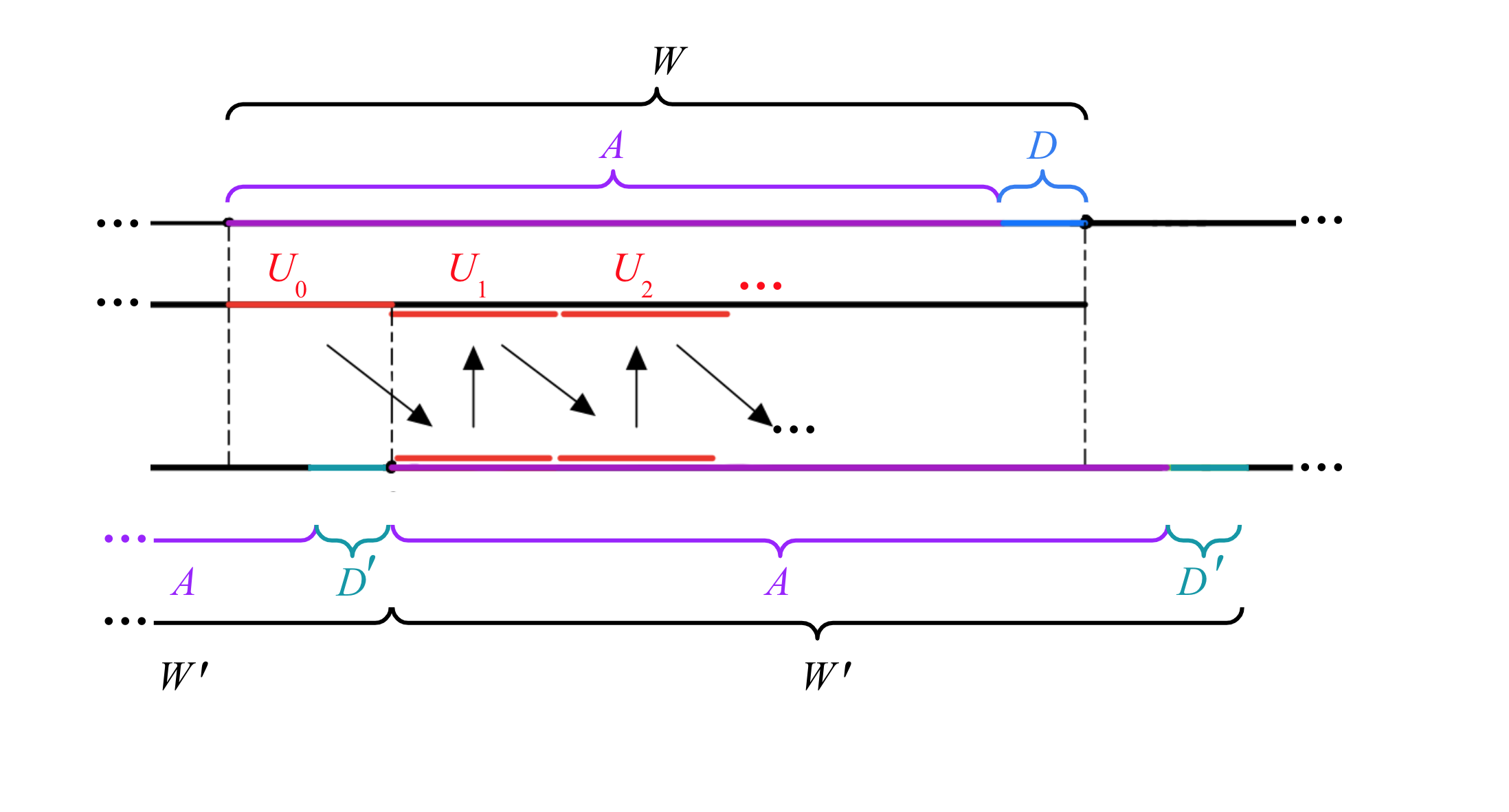}
\caption[]{\label{jump1}
Since $U$ is a prefix in $AD$, it is a prefix in $AD'$. The rightward shift of $AD'$ by $U$ induces a cascade $U=U_0, U_1, U_2, \ldots$ in $AD$, where the $U_i$ concatenate and $U_i\equiv U$ for each $i$.}
\end{figure} 
\end{defn}
\begin{prop}\label{prop:unique jump} Let $W\in F\left(X\right)$ be a cyclically reduced nonperiodic word. Then the maximal ascent in $W$ is uniquely positioned.
\end{prop}
\begin{proof}
    Let $A$ be the maximal ascent in $W$ and suppose without loss of generality that $W=AD$. If $D=1_\mathcal{F}$, then $W=A$ and so $A$ is uniquely positioned in $W$ since 
$W$ is nonperiodic and thus all its cyclic permutations are distinct. If $D\neq 1_\mathcal{F}$, then by Lemma~\ref{lem: ascent is M minus m} $D$ is a descent. Following Remark~\ref{rem:W square}, suppose $W^2=ADAD$ contains a third occurrence of $A$. Since $A$ appears exactly once in each $W$, the third occurrence of $A$ begins in the first factor $W$ and ends in the second one. Let $W'=AD'$ be a cyclic permutation of $W$ starting with $A$. Then $D'$ is a descent with $|D'|=|D|$. Since $W$ is nonperiodic, it has distinct conjugates, and so $D'\not\equiv D$. Moreover, $W$ has a prefix $U_0$ such that $W'=U_0^{-1}WU_0$. The shift $U_0$ induces a cascade $U_0, U_1,U_{2},  \ldots$ in $W=AD$. 
    
    \underline{\textbf{Claim }}:\label{claim: U is longer than D} $|D|<|U_0|<|A|$ and $|U_0|$ is not a divisor of $|W|$. 
\begin{proof}[Proof of Claim] We have $W'=AD'=U_0^{-1}ADU_0$ with $|U_0|<|W|=|A|+|D|$. If $|U_0|>|A|$, then $W'$, and thus $A$, begins in the interior of $D$. Since ascents and descents do not overlap, the ascent $A$ is internal in $D$, and so $A$ is not unique in $W$, which contradicts Lemma~\ref{lem: ascent is M minus m}. If $|U_0|< |D|$, then the ascent $A$ appearing in $W'$ ends in the interior of $D$ which is a contradiction. If $|U_0|=|D|$, then $A$ and $D$ have a common suffix which is impossible.

Since the subwords $U_{\textcolor{black}{i}}$ \textcolor{black}{are concatenated}, if $\dfrac{|W|}{|U_0|}=k\in\mathbb{N}$, then $W=U_0^k$ which is a contradiction.
    \end{proof}
Note that $U_0=A_1D'$ where 
$A_1$ is an ascent. Indeed, $U_0$ is a prefix of $A$ (in $W$) and so each prefix of $A_1$ is $\succ 1_\mathcal{F}$; and $A_1$ is a suffix of $A$ (in $W'$) and so each suffix of $A_1$ is $\succ 1_\mathcal{F}$. 

Consider the cascade induced by $U_0=A_1D'$. Then $W=AD$ is a proper subword of $U_0^n$ for some $n>1$.\textcolor{black}{Indeed, the cascade ensures that the ascent $A$ in $W$ appears as a subword of $U_0^n$ for some $n>1$. However, since $U_0\equiv A_1D'$ where $A_1$ is an ascent and $D'$ is a descent, and $A_1$ does not overlap with $D$ in $W=AD$, the only possibility is for $A_1$ to be long enough so that the last occurrence of $A_1$ in $W=AD$ must begin in $A$, contain $D$, and end in the interior of the second $W$ factor in $W^2$}. 
Let $U_{n}$ be the first term that is not in $W$. Then $U_{n}\equiv A_1D'$ and $W=AD$ overlap.  So $D$ appears as a subword of $U_{0}$. For $D$ to appear as a subword of $U_0$, it is necessary that $D$ appears as an internal subword of $A_1$, since $D$ is not equivalent to $D'$ and $D$ does not overlap with the ascent $A_1$. So we have a new overlap of two subwords equivalent to $U_{1}$. See Figure~\ref{jump3}.
\begin{figure}\centering
\includegraphics[width=.8\textwidth]{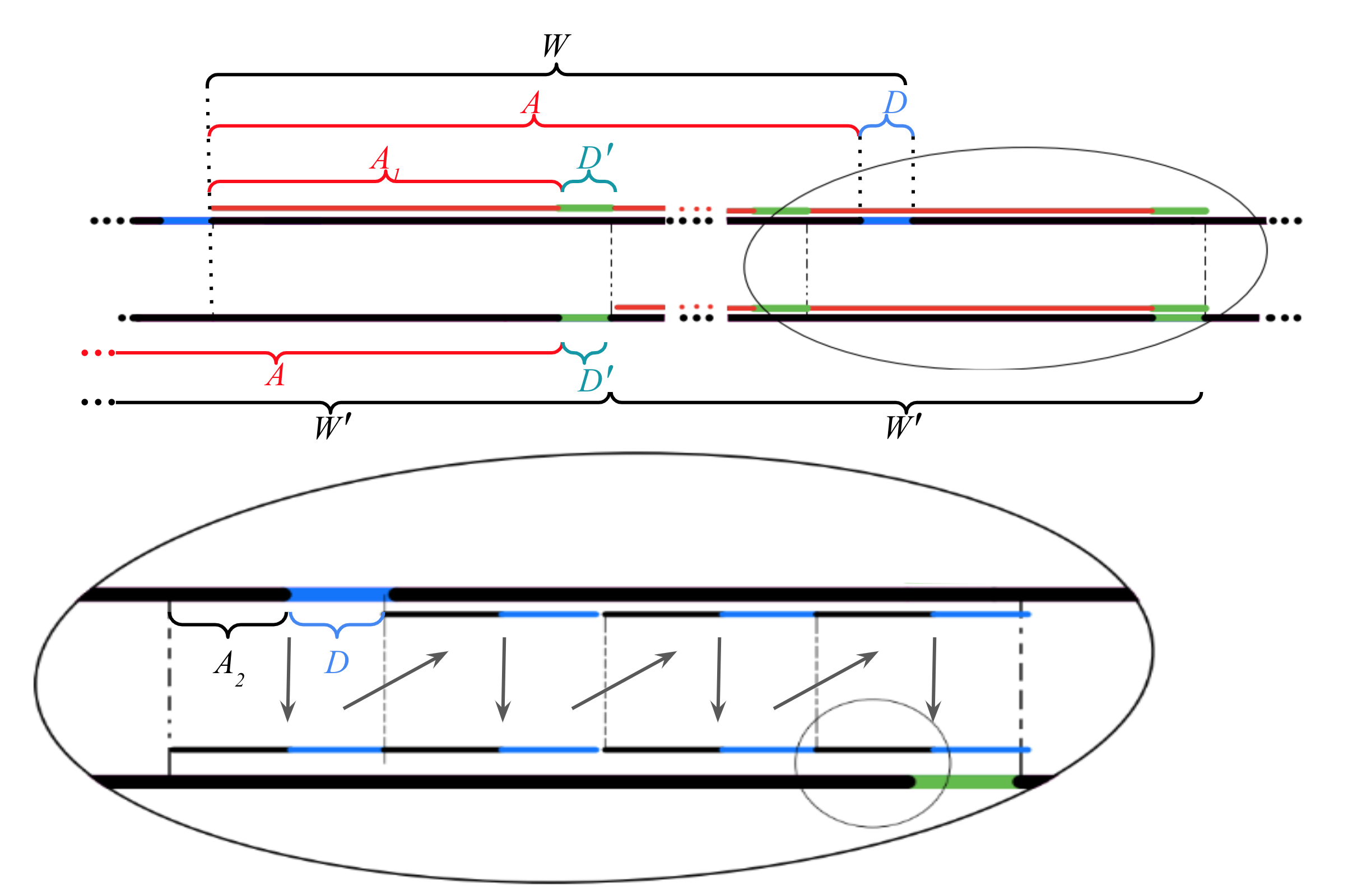}
\caption[]{\label{jump3}
The cascade induced by $U_0=A_1D'$ terminates with $D$ appearing as an internal subword of $A_1$. See the region inside the ellipse. So $A_1$ has a prefix $U_2=A_2D$ which itself is a shift that induces a cascade in $A_1D'$. In this example, the cascade induced by $A_2D$ immediately leads to a contradiction since it forces an overlap between the ascent $A_2$ and the descent $D'$. See the region inside the circle.}
\end{figure} 
The shift $U_2=A_2D$  induces a new cascade that follows the same pattern as above with the difference being that $U_2=A_2D$ is a concatenation of an ascent $A_2$, with $|A_2|<|A_1|<|A|$, and $D$ instead of $D'$. Note that by the above Claim, $\textcolor{black}{|D'| =}\ |D|<|U_2|<|A_1|$. Once again, the cascade of copies of $A_2D$ requires that the copies of $D$ must not coincide with $D'$ and cannot overlap with $A_2$. So $D$ appears as a subword of $A_2$. Thus $A_2$ contains a subword $U_3=A_3D'$ where $|A_3|<|A_2|$ is an ascent. As this process repeats, the shift $U_j$ will be a concatenation of an ascent $A_j$ and the descent $D$ if $j$ is even, and $U_j$ will be a concatenation of an ascent $A_j$ and the descent $D'$ if $j$ is odd. For each cascade, the ascent $A_j$ in the shift $U_j=A_jD$ (or $U_j=A_jD'$) has shorter length. Since $|W|<\infty$ and by induction, this process terminates with either an overlap of an ascent and a descent, or by forcing $D=D'$. Both lead to contradictions.
\end{proof}




\section{Magnus Bi-order}
The following describes an explicit bi-order on free groups due to Magnus~\cite{MKS66}. Let $\mathcal{F}=\mathcal{F}\left(x_1,x_2\right)$ be the free group on generators $x_1, x_2$. Let $\Lambda=\mathbb{Z}[[X_1,X_2]]$ be the ring of formal power series in the non-commuting variables $X_1$ and $X_2$, one for each generator of $\mathcal{F}$.   Define the multiplicative homomorphism $\mu: \mathcal{F}\rightarrow \Lambda$ as:
\[ \mu:\begin{cases} 
      x_i\quad \mapsto 1+X_i \\
      x_i^{-1}\ \mapsto 1-X_i+X_i^{2}-X_i^3+\cdots   
   \end{cases}
\]
For example: 
\begin{align*}
\mu\left(x_1x_2^{-1}\right)&=(1+X_1)(1-X_2+X_2^2-X_2^3+\cdots)\\
&=1+X_1-X_2+O\left(2\right)
\end{align*}
where $O\left(n\right)$
refers to the sum of all terms of order $\geq n$. Then $\mu$ is injective and $\mathcal{F}$ embeds in the group of units $1+O\left(1\right)\subset \Lambda$. Order the elements of $\Lambda$ as follows. First adopt the convention of writing the elements of $\Lambda$ in standard form starting from lower degree terms in an increasing order. Then, order the terms with the same degree lexicographically where $X_1\succ X_2$. Compare two elements of $\Lambda$ according to the coefficients of the first term at which they differ.
For example, $1+X_1+3X_2+O(2)\ \succ\ 1+X_1+X_2+O(2)$ since the first term at which they differ is $X_2$, and the coefficient of $X_2$ in the first element is greater than the coefficient of $X_2$ in the second one. Under this order, $1+O(1)\subset \Lambda$ is a bi-ordered group. 

Define an ordering $\succ$ on $\mathcal{F}$ by:
$$v\ \succ\  w \ \iff\ \mu\left(v\right)\ \succ \ \mu\left(w\right)$$
It is readily verified that $\succ$ is both left and right invariant.
\begin{defn}\label{defn:monotonic} A word $W=y_1\cdots y_n\in \mathcal{F}\left(X\right)$ is \textit{monotonic} if either $y_i\in X$ for each 
 $1\leq i\leq n$, or $y_i\in X^{-1}$ for each 
 $1\leq i\leq n$.
\end{defn}
\begin{thm}\label{thm:main}
Let $X=\left\{x_1, x_2\right\}$ be an alphabet and let\textcolor{black}{ $\mathcal{F}=\mathcal{F}\left(X\right)$ be the free group on $X$ equipped with a bi-order $\prec$. Let $W\in \mathcal{F}$} be a cyclically reduced nonperiodic word of length $>1$. Then $W$ has a cyclic permutation $W'=AD$ where:
\begin{enumerate}
\item $A$ is the uniquely positioned maximal ascent in $W$.
\item If $D$ is not uniquely positioned, then it appears as an internal subword of $A$.
\item Using the Magnus ordering on $\mathcal{F}$, we have $D=1_\mathcal{F}$ if and only if $W$ is monotonic.
\end{enumerate}
\end{thm}
\begin{proof}
Let $A$ be the maximal ascent in $W$. By Proposition~\ref{prop:unique jump}, $A$ is uniquely positioned. 

Let $W'=AD$\textcolor{black}{where $D$ is not uniquely positioned in $W$. If $D\neq 1_\mathcal{F}$, then by} Lemma~\ref{lem: ascent is M minus m}, $D$ is a descent. Suppose $D'\equiv D$ is a subword in $AD$. By Remark~\ref{rem:no overlap}, $D'$ has no overlap with $A$, and so it appears as an internal subword of $A$.\textcolor{black}{Moreover, if $D=1_\mathcal{F}$, then $W=A$ and $D$ is the empty word between any concatenated subwords of $W$. Thus $D$ appears as an internal, albeit trivial, subword of $W$. Note that by assumption, $|W|>1$ and so $W$ has at least two subwords.}

Choose the Magnus bi-ordering corresponding to $x_1\succ x_2\succ 1_\mathcal{F}$. Then any nonempty word in $X$ is $\succ 1_\mathcal{F}$. If $W$ is monotonic, then so is each cyclic permutation of $W$. Suppose without loss of generality that $W'=y_1\cdots y_n$ with $y_i\in X$, for $1\leq i\leq n$.  By the maximality of $A$, if $D\neq 1_\mathcal{F}$, then it is a descent, which is impossible since all monotonic words are $\succ 1_\mathcal{F}$.

Suppose $D=1_\mathcal{F}$. Then $W'$ is the maximal ascent. Suppose $W'$ is not monotonic. Then $W=Ux_1^{-1}V$ for some words $U, V\in \mathcal{F}$. The case $W=Ux_2^{-1}V$ is similar. Let $$\mu(U)=1+M_1X_1+M_2X_2+O(2)\quad \text{and}\quad \mu(V)=1+N_1X_1+N_2X_2+O(2)$$ 
where $M_i, N_i \in \mathbb{Z}$.
Then $$\mu\left(VU\right)=1+(M_1+N_1)X_1+(M_2+N_2)X_2+O(2)$$
Moreover, we have 
\begin{align*}
\mu\left(W'\right)&=(1+M_1X_1+M_2X_2+O(2))(1-X_1+O(2))(1+N_1X_1+N_2X_2+O(2))\\
&=1+(M_1+N_1-1)X_1+(M_2+N_2)X_2+O(2)
\end{align*}
Hence $W'\prec VU$. But $VU$ is a prefix of the cyclic permutation $x_1U^{-1}W'Ux_1^{-1}=VUx_1^{-1}$ contradicting the maximality of $W'$.
\end{proof}

\bibliographystyle{alpha}
\bibliography{brahim.bib}

\end{document}